\documentclass[11pt, reqno, a4paper]{article}
\usepackage{amsfonts}
\usepackage{mathrsfs}
\usepackage{amssymb}
\usepackage{amsmath}
\newtheorem{theorem}{Theorem}[section]
\newtheorem{lemma}{Lemma}[section]
\newtheorem{corollary}{Corollary}[section]
\newtheorem{definition}{Definition}[section]

\begin{document}

\title{Approximation by  Boolean Sums of Jackson Operators on the Sphere
\thanks{The research was supported by the National
Natural Science Foundation of China (No. 60873206), the Natural
Foundation of Zhejiang Province of China (No. Y7080235) and the
Innovation Foundation of Post-Graduates of Zhejiang Province of
China (No. YK2008066)}}

\author{ Yuguang Wang
\and
 Feilong Cao\thanks{Corresponding author.
 E-mail: \tt feilongcao@gmail.com}}

  \date{}
\maketitle

\begin{center}
\footnotesize

 Institute of Metrology and Computational Science,
 China Jiliang University,

 Hangzhou 310018, Zhejiang Province,  P. R. China.

\begin{abstract}
This paper concerns the approximation by the Boolean sums of Jackson operators $\oplus^rJ_{k,s}(f)$ on the unit sphere $\mathbb S^{n-1}$ of $\mathbb{R}^{n}$.
We prove the following the direct and inverse theorem for $\oplus^rJ_{k,s}(f)$: there are constants $C_1$ and $C_2$ such that
\begin{equation*}
C_1\|\oplus^rJ_{k,s}f-f\|_p \leq \omega^{2r}(f,k^{-1})_p \leq C_2
\max_{v\geq k}\|\oplus^rJ_{k,s}f-f\|_p
\end{equation*}
for any positive integer $k$ and any $p$th Lebesgue integrable functions $f$ defined on $\mathbb S^{n-1}$, where $\omega^{2r}(f,t)_p$ is the modulus of
smoothness of degree $2r$ of $f$. We also prove that the saturation order for $\oplus^rJ_{k,s}$ is $k^{-2r}$.


\end{abstract}
\end{center}

\section{Introduction}
In past decades, many mathematicians have dedicated to establish the
Jackson and Bernstein-type theorems on the sphere. Lizorkin and
Nikol'ski\v{\i} \cite{Lizorkin1983}  constructed   Boolean sums of
Jackson operators
$\oplus^rJ_{k,s}$ (which will be defined in the next section)
 for proving the direct and
inverse theorems on a  special Banach space
$H_p^r(\mathbb S^{n-1})$.
Later, by using a modulus of smoothness as metric,
 Lizorkin and Nikol'ski\v{\i} \cite{Lizorkin1986}
 obtained the direct estimate for Jackson operators (i.e., the $1$-th
Boolean sums of Jackson operators) approximating continuous
function defined on the unit sphere $\mathbb S^{n-1}$.
In 1991,  Li and Yang  \cite{Li1991}
 used the equivalent relation between the $K$-functional and
the smoothness (see \cite{Butzer1971}) and established an inverse
inequality of weak type, which was also called Steckin-Marchaud type
inequality.
 Recently, Ditzian
\cite{Ditzian2004} proved an equivalence relation between
$K$-functional and modulus of smoothness  of high order, which will
provide a tool and allow us to make a proof for direct and inverse
theorems for the  Boolean sums of Jackson operators.

Actually, after improving a
Steckin-Marchaud type inequality, we will prove the direct and inverse
theorem of approximation for arbitrary $r$-th Boolean sums of
Jackson operators $\oplus^rJ_{k,s}$ approximating
$p$-th Lebesgue integrable function on the sphere. Particularly,
 a converse inequality of strong type ({see \cite{Ditzian1994})
 for   $\oplus^rJ_{k,s}$
will be established.
  Moreover, we will use the method of
multipliers and obtain the saturation order of $\oplus^rJ_{k,s}$.

 \section{Definitions and Auxiliary Notations}
Let $\mathbb{S}^{n-1}$ be the unit sphere in Euclidean space
$\mathbb{R}^n$. We denote by the letters $C$ and $C_i$ positive
constants, where $i$ is either positive integers or variables on
which $C$ depends only. Their values may be different at different
occurrences, even within the same formula. We shall denote by $x$
and $y$ the points of $\mathbb {S}^{n-1}$. The notation $a \approx
b$ means that there exists a positive constant $C$ such that
$C^{-1}b\leq a \leq Cb$.

We now introduce some concepts and properties of sphere (see also
\cite{Muller1966}, \cite{Wang2000}).
 The volume of
$\mathbb{S}^{n-1}$
is
$$
\left|\mathbb{S}^{n-1}\right|:=\int_{\mathbb{S}^{n-1}}d\omega=\displaystyle\frac{2\pi^{\frac{n}2}}{\Gamma(\frac{n}2)}.
$$
Denote by $L^p(\mathbb{S}^{n-1})$ the space of $p$-th integrable
functions on $\mathbb{S}^{n-1}$ endowed with the norms
$$
\| f\|_{\infty}:=\|f\|
_{L^\infty(\mathbb{S}^{n-1})}:=\mbox{ess}\!\!\sup_{x\in
\mathbb{S}^{n-1}}|f(x)|
$$
and
$$
\|f\|_p:=\|f\|_{L^p(\mathbb{S}^{n-1})}:=\left\{\int_{\mathbb{S}^{n-1}}|f(x)|^pd\omega(x)\right\}^{1/p}<\infty,
\quad
1\leq p<\infty.
$$

For $f\in L^1(\mathbb{S}^{n-1})$,  the translation operator is
defined by (see for instance \cite{Wang2000})
$$
S_\theta(f)(x):=\frac1{|\mathbb{S}^{n-2}|\sin^{n-2}\theta}\int_{x\cdot
y=\cos \theta}f(y)d\omega'(y),\quad 0<\theta<\pi
$$
where $d\omega'(y)$ denotes the elementary surface piece on
$\mathbb{S}^{n-2}$.\\
We set \[S_\theta^{(0)}(f):=f, \;
S_\theta^{(j)}(f):=S_\theta S_\theta^{(j)}(f),\quad j=1,2,3,\dots\]
and introduce the spherical differences (see \cite{Lizorkin1983})
\[\Delta_\theta^1(f):=\Delta_\theta(f):=S_\theta(f)-f\]
and
\begin{eqnarray*}
\Delta_\theta^r(f):=\Delta_\theta\Delta_\theta^{r-1}(f)=\sum_{j=0}^r
(-1)^{r+j} {r \choose j}
S_\theta^{(j)}(f)=\left(S_\theta-I\right)^r f,\quad
r=2,3,\dots
\end{eqnarray*}
where $I$  is the identity operator on $L^p(\mathbb S^{n-1})$.
Then the modulus of smoothness of degree $2r$ of $f\in L^p(\mathbb{S}^{n-1})$
is defined by (see \cite{Rustamov1994})
\begin{eqnarray*}
\omega^{2r}(f,t):=\sup_{0<\theta\leq t}\|\Delta_\theta^rf\|_p,\quad
0<t<\pi,\; r=1,2,\dots.
\end{eqnarray*}

We denote by $\widetilde{\Delta}$ the Laplace-Beltrami operator
$$
\widetilde{\Delta} f:=\sum_{i=1}^n\frac{\partial^2g(x)}{\partial
x_i^2}\bigg|_{|x|=1},\quad
g(x):=f\left(\frac{x}{|x|}\right),\quad\left(f,\;\widetilde{\Delta}
f\in L^p(\mathbb{S}^{n-1})\right)
$$
by which  we  introduce $K$-functional of degree $2r$
on
$\mathbb{S}^{n-1}$ as
\begin{eqnarray*}
K_{2r}(f,\widetilde{\Delta},t^{2r})_p:=\inf\left\{\|f-g\|_{p}+t^{2r}\|\widetilde{\Delta}^r
g\|_p:\ g,\; \widetilde{\Delta}^r g\in
L^p(\mathbb{S}^{n-1})\right\}.
\end{eqnarray*}
For the modulus of smoothness and $K$-functional, Ditzian
\cite{Ditzian2004} has obtained the following equivalence relation
\begin{eqnarray}\label{sec2eq5}
\omega^{2r}(f,t)_p \approx K_{2r}(f,\widetilde{\Delta},t^{2r})_p.
\end{eqnarray}
Clearly, for any $\lambda>0$, we have, using (\ref{sec2eq5}),
\begin{eqnarray}\label{sec2eq6}
\omega^{2r}(f,\lambda t)_p\leq C_1
K_{2r}(f,\widetilde{\Delta},(\lambda t)^{2r})_p & \leq & C_1
\max\{1,\lambda^{2r}\} K_{2r}(f,\widetilde{\Delta},t^{2r})_p\nonumber\\
& \leq & C_2 \max\{1,\lambda^{2r}\}\omega^{2r}(f,t)_p
\end{eqnarray}
where $C_1$ and $C_2$ are independent of $t$ and $\lambda$.

 Spherical polynomials of order $k$ on $\mathbb{S}^{n-1}$ is defined
 by (see \cite{Lizorkin1983})
\begin{eqnarray*}
P_k(x):=\sum_{j=0}^k H_j(x)
\end{eqnarray*}
where $H_j(x)$ is the spherical harmonic of order $j$, the trace of
some homogeneous polynomials on $\mathbb{R}^{n}$
\[Q_j(x):=\sum_{i_1+i_2+\dots+i_n=m}\!\!\!\!\!\!\!\!x_1^{i_1}\cdots x_n^{i_n}\]
of order $j$, where $x_i (i=1,2,\dots,n)$ is the $i$-th coordinate
of $x\in \mathbb{R}^{n}$. We denote by $\Pi_k^{n}$ the collection of
all spherical polynomials on $\mathbb{S}^{n-1}$ of degree no more
than $k$.

The known Jackson operators are defined by (see \cite{Lizorkin1983})
\begin{eqnarray*}
J_{k,s}(f)(x):= \frac{1}{|\mathbb{S}^{n-2}|} \int_{\mathbb{S}^{n-1}}
f(y)\mathscr{D}_{k,s}(\arccos x\cdot y) d\omega(y),\quad f\in
L^p(\mathbb{S}^{n-1})
\end{eqnarray*}
where $k$ and $s$ are positive integers, $d\omega(y)$ is the
elementary surface piece of $\mathbb{S}^{n-1}$,
$$\mathscr{D}_{k,s}(\theta):=A^{-1}_{k,s}\left({\displaystyle\frac{\sin\frac{k\theta}{2}}{\sin\frac{\theta}{2}}}\right)^{2s}$$
is the classical  Jackson kernel where $A^{-1}_{k,s}$ is
a constant connected with $k$ and $s$ such that
\[
\int_0^\pi
\mathscr{D}_{k,s}(\theta)\sin^{2\lambda}\theta\:d\theta=1,
\quad
\lambda=\frac{n-2}{2}.
\]
We observe that
\begin{eqnarray}\label{sec2eq7}
J_{k,s}(f)(x) & = & \frac{1}{|\mathbb{S}^{n-2}|}
\int_{\mathbb{S}^{n-1}} f(y)\mathscr{D}_{k,s}(\arccos x\cdot
y) d\omega(y)\nonumber\vspace{0.1 cm}\\
& = &
\int_0^\pi\left(\frac{1}{|\mathbb{S}^{n-2}|\sin^{2\lambda}\theta}
\int_{x\cdot y=\cos\theta}
f(y)d\omega'(y)\right)\mathscr{D}_{k,s}(\theta)\sin^{2\lambda}\theta\:d\theta\nonumber\\
& = & \int_0^\pi
S_\theta(f)(x)\mathscr{D}_{k,s}(\theta)\sin^{2\lambda}\theta\:d\theta.
\end{eqnarray}

We  introduce the definition
of $r$-th  Boolean sums of Jackson operator  as follows (see
\cite{Lizorkin1983}).

\begin{definition}\label{sec2def1}
The r-th $(r\geq1)$   Boolean sum of Jackson operator of degree $k$
($k\geq1$) on $\mathbb S^{n-1}$ is defined by
\begin{eqnarray}\label{sec2eq8}
\oplus^rJ_{k,s}(f):=\left(-(I-J_{k,s})^r+I\right)(f),\quad f\in
L^p(\mathbb{S}^{n-1}),
\end{eqnarray}
where $s$ is a positive integer.
\end{definition}

It is clear that
\begin{eqnarray}\label{sec2eq10}
\oplus^rJ_{k,s}(f)=-\sum_{i=1}^r (-1)^i {r \choose i} (J_{k,s})^i(f).
\end{eqnarray}

 We now make a brief introduction
of projection operators $Y_{j}(\cdot)$ by ultraspherical
(Gegenbauer) polynomials
$\{G_{j}^{\lambda}\}_{j=1}^{\infty}$
 $(\lambda={\frac{ n-2}{2}})$ for discussion
of saturation property of $\oplus^rJ_{k,s}$.

Ultraspherical polynomials $\{G_j^{\lambda}\}_{j=1}^{\infty}$ are defined in terms of the generating function (see \cite{Stein1971}):
\[\frac{1}{(1-2 t r+r^2)^\lambda}=\sum_{j=0}^{\infty}G_j^{\lambda}(t)r^j\]
where $|r|<1$, $|t|\leq 1$. For any $\lambda>0$, we have (see \cite{Stein1971})
\begin{equation}\label{eq24}
G_1^{\lambda}(t)=2\lambda t
\end{equation} and
\begin{equation}\label{eq27}
\displaystyle\frac{d}{dt}G_j^{\lambda}(t)= 2\lambda
G_{j-1}^{\lambda+1}(t).
\end{equation}
When $\lambda=\frac{n-2}{2}$ (see \cite{Wang2000}),
\[G_j^{\lambda}(t)=\frac{\Gamma(j+2\lambda)}{\Gamma(j+1)\Gamma(2\lambda)}P_j^n(t),\quad
j=0,1,2,\dots\] where $P_j^n(t)$ is the Legendre polynomial of
degree $j$ (see \cite{Muller1966}).
Particularly,
\[G_j^{\lambda}(1)=\frac{\Gamma(j+2\lambda)}{\Gamma(j+1)\Gamma(2\lambda)}P_j^n(1)=\frac{\Gamma(j+2\lambda)}{\Gamma(j+1)\Gamma(2\lambda)},\quad
j=0,1,2,\dots,\]therefore,
\begin{equation}\label{eq28}
P_j^n(t)=\frac{G_j^{\lambda}(t)}{G_j^{\lambda}(1)}.
\end{equation}
Besides, for any $j=0,1,2,\ldots$, and $|t|\leq 1$,
$|P_j^n(t)|\leq1$ (see \cite{Muller1966}).

The projection operators are defined by
\begin{eqnarray*}
Y_j(f)(x):=\frac{\Gamma(\lambda)(n+\lambda)}{2\pi^{\frac{n}{2}}}\int_{\mathbb{S}^{n-1}}G_j^{\lambda}(x\cdot
y)f(y)\:d\omega(y).
\end{eqnarray*}
It follows from (\ref{eq24}), (\ref{eq27}) and (\ref{eq28}) that
\begin{eqnarray}\label{eq25}
\lim\limits_{t\rightarrow 1-}
{\displaystyle\frac{1-P_j^n(t)}{1-P_1^n(t)}=\frac{j(j+2\lambda)}{2\lambda+1}},\quad
j=0,1,2,\dots.
\end{eqnarray}

Finally, we introduce the definition of saturation for operators
(see \cite{Berens1968}).
\begin{definition}
 Let $\varphi(\rho)$ be a positive function with respect to $\rho$,
$0<\rho<\infty$, tending monotonely to zero as
$\rho\rightarrow\infty$. For a sequence of
operators $\{I_\rho\}_{\rho>0}$ if there exists $\mathcal {K}\subseteqq
L^p(\mathbb S^{n-1})$ such that
\vspace{0.1 cm}\\\parbox{6
cm}{\begin{tabular}{ll}
(i) & If\; $\| I_\rho(f)-f\|_{p} = o(\varphi(\rho))$, then $I_\rho(f)=f$; \\
(ii) & $\| I_\rho(f)-f\|_{p} = O(\varphi(\rho))$ if and only if
$f\in\mathcal {K}$;
\end{tabular}}\parbox{1 cm}{}\vspace{0.2 cm}\\
then $I_\rho$ is said to be saturated on $L^p(\mathbb S^{n-1})$ with
order $O(\varphi(\rho))$ and $\mathcal {K}$ is called its saturation
class.
\end{definition}

\section{Some Lemmas}
In this section, we show some lemmas as the preparation for the
proof of the main results.
\begin{lemma}\label{lm31}
For any $f\in L^p(\mathbb{S}^{n-1})$ and any positive integers $k$, $r$, $u$,
$s$, we have,
\begin{tabular}{ll}
 (i)    & $\oplus^uJ_{k,s}(f)$ is a spherical polynomial of degree no more than $s(k-1)$, which\\
   &implies $\widetilde{\Delta}^{r}(\oplus^uJ_{k,s}(f))\in L^p(\mathbb{S}^{n-1})$;\vspace{0.1 cm}\\
 (ii)   & $\|\oplus^uJ_{k,s}(f)\|_p\leq 2^u\|f\|_p$;\vspace{0.1 cm}\\
 (iii)  & $\|\widetilde{\Delta}^{r}(\oplus^uJ_{k,s}(f))\|_{p}\leq Ck^{2r}\|
 f\|_{p}$;\vspace{0.1 cm}\\
 (iv)   & If \;$\widetilde{\Delta}^{r}g\in
 L^p(\mathbb{S}^{n-1})$,\; then $\|\widetilde{\Delta}^{r}(\oplus^uJ_{k,s}(g))\|_{p}\leq 2^u\| \widetilde{\Delta}^{r}g\|_{p}$.
\end{tabular}
\end{lemma}

{\bf Proof.} (i) Since $\mathscr{D}_{k,s}(\theta)$ is an even
trigonometric polynomial of degree no more than $k(s-1)$, then
$J_{k,s}(f)$ is a spherical polynomial of degree no more than
$s(k-1)$.
Thus we can prove by induction that
$\oplus^uJ_{k,s}(f)$
is a
spherical polynomial of degree no more than $s(k-1)$.

(ii) Using the contraction of translation operator  that (see for
instance \cite{Wang2000})
\begin{eqnarray}\label{sec3eq1}
\|S_\theta(f)\|_p\leq \|f\|_p,\quad0<\theta<\pi
\end{eqnarray}
as well as (\ref{sec2eq7}), we have
\begin{eqnarray*}
\|J_{k,s}(f)\|_p &\leq & \int_0^\pi
\|S_\theta(f)\|_p\mathscr{D}_{k,s}(\theta)\sin^{2\lambda}\theta\:d\theta\leq
\|f\|_p.
\end{eqnarray*}
Using (\ref{sec2eq10}), we get
\begin{eqnarray*}
\|\oplus^uJ_{k,s}(f)\|_p &\leq & 2^u \:\|f\|_p.
\end{eqnarray*}

(iii)  Using Bernstein-inequality on the sphere (see
\cite{Pawelke1972}), that is,
\[\|\widetilde{\Delta}P_k\|_{p}\leq Ck^{2}\| P_k\|_{p},
\quad
P_k\in \Pi^n_k,
\]
 we may
easily obtain by induction that
\begin{eqnarray*}
\|\widetilde{\Delta}^{r}P_k\|_{p}\leq C^{r}k^{2r}\| P_k\|_{p}.
\end{eqnarray*}
Noting that (i) implies $\oplus^uJ_{k,s}(f)\in\Pi^{d}_{s(k-1)}$, we thus have
\begin{eqnarray*}
\|\widetilde{\Delta}^{r}(\oplus^uJ_{k,s}(f))\|_{p} &\leq& C
k^{2r}\|\oplus^uJ_{k,s}(f)\|_p
 =  C_r k^{2r}\| f \|_{p},
\end{eqnarray*}
where $C$ and $C_r$ are independent of $f$ and $k$.

(iv) The fact that
$\widetilde{\Delta}S_{\theta}(g)=S_{\theta}\widetilde{\Delta}(g)$
(see \cite{Lizorkin1983}) implies
\begin{eqnarray}\label{sec3eq3}
\widetilde{\Delta}^r( {\oplus^u}J_{k,s}(g))=
{\oplus^u}J_{k,s}(\widetilde{\Delta}^r g).
\end{eqnarray}
We thus use (ii) and find
$$
\|\widetilde{\Delta}^{r}(\oplus^uJ_{k,s}(g))\|_{p}=\|\oplus^uJ_{k,s}(\widetilde{\Delta}^{r}g)\|_{p}\leq
2^u \|\widetilde{\Delta}^{r}g\|_{p}.
$$
This completes the
proof of Lemma 3.1.
\quad$\Box$

\begin{lemma}\label{lm32}
For $\beta\geq-1$, $2s\geq\beta+n-2$, ${0<\gamma\leq\pi}$, and
$n\geq3$,  we have
\begin{eqnarray}\label{eq34}
\int_0^\gamma \theta^\beta
\mathscr{D}_{k,s}(\theta)\sin^{2\lambda}\theta\:d\theta\approx
k^{-\beta},
\end{eqnarray}
where ${\lambda=\frac{n-2}{2}}$, and $s$, $n$, $k$ are positive integers.
\end{lemma}
{\bfseries Proof.} A simple computation gives
\begin{eqnarray*}
\int_{0}^{\gamma}\theta^{\beta}\mathscr{D}_{k,s}(\theta)\sin^{2\lambda}\theta
d\theta =
\frac{\displaystyle\int_{0}^{\gamma}\theta^{\beta}\left(\frac{\sin\frac{k\theta}{2}}{\sin\frac{\theta}{2}}\right)^{2s}\sin^{2\lambda}\theta
d\theta}{\displaystyle\int_{0}^{\gamma}\left(\frac{\sin\frac{k\theta}{2}}{\sin\frac{\theta}{2}}\right)^{2s}\sin^{2\lambda}\theta
d\theta} &\approx& k^{-\beta}.
\end{eqnarray*}

The proof of Lemma~\ref{lm32} is completed.
\quad$\Box$

The following lemma is an improved version of Lemma 2.1 of
\cite{Wickeren1986}, which is useful for the proof of Bernstein-type
inequality.

\begin{lemma}\label{sec3lm3}
Assuming that
$\left\{\sigma_v\right\}_{v=1}^\infty,\left\{\tau_v\right\}_{v=1}^\infty$
are nonnegative sequences with $\sigma_{1}=0$. For $a > 1$ and
positive integer $k$, if the inequalities
\[\sigma_{k}\leq\left(\frac{av}{k}\right)^{p}\sigma_{v}+\tau_{v},\quad
v=1,2,\dots,k\] hold, then there exists $0<b_{m}<1$,
$b_{m}\rightarrow 1$, as $m \rightarrow \infty$, such that
\[\sigma_{k}\leq C_m
k^{-b_{m}p}\sum^{k}_{v=1}v^{b_{m}p-1}\tau_{v}.\]
\end{lemma}

{\bf Proof.} First we take $m>1$, such that
$\displaystyle\ln\frac{m}{a}>1$. Then there exists positive integer
$N$, such that $m^{N}\leq k< m^{N+1}$. Take $s_{v}\geq0$ such that
$km^{-v-1}< s_{v} \leq km^{-v}$, $v=0,1,\dots,N$ as well as
\[\tau_{s_v}\leq \tau_{j}\quad(km^{-v-1}<j\leq km^{-v}).\]
Set $s_{N+1}=1$, $b_{m}=\displaystyle\frac{\ln\frac{m}{a}}{\ln m}<
1$. Clearly $b_{m}\rightarrow 1$, as $m\rightarrow \infty$. Then
\begin{eqnarray*}
\sigma_{k}&\leq &
\left(\frac{a s_0}{k}\right)^{p}\sigma_{s_0}+\tau_{s_{0}}\\
&\leq &
k^{-p}\sum^{N}_{v=0}(a^{v+1}s_v)^{p}\left(\sigma_{s_{v}}-\left(\frac{a s_{v+1}}{s_{v}}\right)^{p}\sigma_{s_{v+1}}\right)+\tau_{s_{0}}\\
&\leq &
m^{p}\sum^{N}_{v=0}\left(\displaystyle\frac{m}{a}\right)^{-p(v+1)}\tau_{s_{v+1}}+\tau_{s_{0}}\\
&\leq & m^{p}\sum^{N+1}_{v=0}m^{-\left(\frac{\ln\frac{m}{a}}{\ln
m}\right)pv}\tau_{s_v}\\
&\leq & C_m k^{-b_{m}p}\left(\sum^{N+1}_{v=0}\sum_{km^{-v-1}< j
\leq km^{-v}}j^{b_{m}p-1}\tau_{j}+\tau_1\right)\\
&= & C_m k^{-b_{m}p}\sum^{k}_{j=1}j^{b_{m}p-1}\tau_{j}.
\end{eqnarray*}

This finishes the proof of Lemma~\ref{sec3lm3}.\quad$\Box$

The following lemma gives the description of $\oplus^rJ_{k,s}$ by
multipliers, which was proved in \cite{Lizorkin1983}.

\begin{lemma}\label{lm34}
For $f\in L^{p}(\mathbb S^{n-1})$, there holds
\begin{eqnarray*}
\oplus^rJ_{k,s}(f) & = & \sum_{j=0}^{\infty}{^r\xi_{k,s}(j)}Y_j(f)
\end{eqnarray*}
where \[^r\xi_{k,s}(j)=1-\left(\int_0^{\gamma}
\mathscr{D}_{k,s}(\theta)\left(1-P_j^n(\cos\theta)
\right)\sin^{2\lambda}\theta\:d\theta\right)^r,\quad j=0,1,2,\dots
\]and the convergence of the series is meant in a weak sense.
\end{lemma}

The final lemma is useful for determining the saturation order.
It can be deduced by the methods in \cite{Berens1968} and
\cite{Butzer1972}.

\begin{lemma}\label{lm38}
Suppose that $\left\{I_\rho\right\}_{\rho>0}$ is a sequence of
operators on $L^p(\mathbb{S}^{n-1})$, and there exists series
$\left\{\lambda_{\rho}(j)\right\}_{j=1}^{\infty}$ with respect to
$\rho$, such that
\[I_\rho(f)(x)=\sum_{j=0}^{\infty}\lambda_\rho(j)Y_j(f)(x)\]for every $f\in L^p(\mathbb{S}^{n-1})$. If for any $j=0,1,2,\dots$,
 there $\mathit{exists}$ $\varphi(\rho)\rightarrow
0\!+(\rho\rightarrow\rho_0)$ such that
\[\lim_{\rho\rightarrow\rho_0}\frac{1-\lambda_\rho(j)}{\varphi(\rho)}=\tau_j\neq0,\] then
$\left\{I_\rho\right\}_{\rho>0}$ is saturated on
$L^p(\mathbb{S}^{n-1})$ with the order $O(\varphi(\rho))$ and the
collection of all constants is the invariant class for
$\left\{I_\rho\right\}_{\rho>0}$ on $L^p(\mathbb{S}^{n-1})$.$\quad
\Box$
\end{lemma}

\section{Main Results and Their Proof}

In this section, we shall state and prove the main results, that is,
the lower and upper bounds as well as the saturation order for
Boolean sums of Jackson operators on $L^p(\mathbb{S}^{n-1})$.

\begin{theorem}\label{th41}
 Let $2s\geq n$,
 and let $\{\oplus^rJ_{k,s}\}_{k=1}^{\infty}$ be the sequence of   Boolean sums of Jackson operators
defined above. Then for any positive integers $k$ and $r$ as well as
sufficiently smoothing $g\in L^{p}(\mathbb S^{n-1}), 1\leq p\leq \infty$
 such that
$\widetilde{\Delta}^{r}g\in
 L^{p}(\mathbb S^{n-1})$, we have
\begin{eqnarray}\label{th41eq18}
\|\oplus^rJ_{k,s}(g)-g\|_{p}\leq
C_1k^{-2r}\|\widetilde{\Delta}^{r}g\|_{p},
\end{eqnarray}
 therefore, for any $f\in L^{p}(\mathbb S^{n-1})$, we have
\begin{eqnarray}\label{th41eq1}
\|\oplus^rJ_{k,s}(f)-f\|_{p}\leq C_2\omega^{2r}(f,k^{-1})_{p},
\end{eqnarray}
where $C_1$ and $C_2$ are constants independent of $f$ and $k$.
\end{theorem}

{\bf Proof.} By Definition~\ref{sec2def1}, we have
$$\oplus^rJ_{k,s}(g)(x)-g(x)=-(I-J_{k,s})^{r}(g)(x).$$
Now we prove (\ref{th41eq18}) by induction.
 For $r=1$,
$$S_\theta(g)(x)-g(x)=\int_{0}^{\theta}\sin^{-2\lambda}\tau\int_{0}^{\tau}\sin^{2\lambda}uS_{u}(\widetilde{\Delta}g)(x)dud\tau$$ (see \cite{Pawelke1972})
implies (explained below)
\begin{eqnarray}\label{th41eq2}
& & \| J_{k,s}(g)-g\|_{p}\;=\;
\left\|\int_{0}^{\pi}D_{k,s}(\theta)\left(S_\theta(g)(\cdot)-g(\cdot)\right)\sin^{2\lambda}\theta
d\theta\right\|_{p}\nonumber\\
&\leq&\int_{0}^{\pi}D_{k,s}(\theta)\sin^{2\lambda}\theta\int_{0}^{\theta}\sin^{-2\lambda}\tau\int_{0}^{\tau}\sin^{2\lambda}u\left\|
S_{u}(\widetilde{\Delta}g)\right\|_{p}dud\tau d\theta\nonumber\\
& \leq &
\sup_{\theta>0}\left\{\theta^{-2}\int_{0}^{\theta}\sin^{-2\lambda}\tau\int_{0}^{\tau}\sin^{-2\lambda}u
du\:d\tau\right\}\left(\int_{0}^{\pi}\theta^{2}D_{k,s}(\theta)\sin^{2\lambda}\theta
d\theta\right)\|\widetilde{\Delta}g\|_{p}\nonumber\\
& \leq & Ck^{-2}\|\widetilde{\Delta}g\|_{p},
\end{eqnarray}
where the Minkowski inequality is used in the first inequality, the
second one by (\ref{sec3eq1}) and the third one is deduced from
Lemma~\ref{lm32}.

 Assume that for any fixed positive integer $u$,
$$\|\oplus^uJ_{k,s}(g)-g\|_{p} \leq C
k^{-2u}\|\widetilde{\Delta}^ug\|_{p}.$$ Then
\begin{eqnarray*}
\|\oplus^{u+1}J_{k,s}(g)-g\|_{p} & = &
\left\|(J_{k,s}-I)(\oplus^uJ_{k,s}(g)-g)\right\|_{p}\leq
Ck^{-2}\|\widetilde{\Delta}(\oplus^uJ_{k,s}(g)-g)\|_{p}\\
& = &
Ck^{-2}\|\oplus^uJ_{k,s}(\widetilde{\Delta}g)-\widetilde{\Delta}g\|_{p}
\leq Ck^{-2u-2}\|\widetilde{\Delta}^{u+1}g\|_{p},
\end{eqnarray*}
where the first inequality is by (\ref{th41eq2}), the second
one by (\ref{sec3eq3}), the last by induction assumption.
Therefore, (\ref{th41eq18}) holds.

 Using (\ref{sec2eq5}) and
noticing that $\oplus^uJ_{k,s}$ is a linear operator, we obtain
(\ref{th41eq1}). This completes the proof of the theorem.\quad$\Box$

Next, we establish an inverse inequality of strong type
 for $\oplus^rJ_{k,s}$
on $L^p(\mathbb{S}^{n-1})$.

\begin{theorem}\label{th42}
For positive $r\geq 1$  and
$f\in L^p(\mathbb S^{n-1})$, $1\leq p\leq \infty$,
there exists a constant $C$ independent of $f$
and $k$ such that
\begin{eqnarray}\label{sec4eq21}
\omega^{2r}(f,k^{-1})_p\leq C \max_{v\geq
k}\|\oplus^rJ_{v,s}(f)-f\|_{p}.
\end{eqnarray}
\end{theorem}

{\bf Proof.} We first establish a Steckin-Marchaud type inequality,
that is, for $f\in L^{p}(\mathbb{S}^{n-1})$,
$$
\omega^{2r}(f,k^{-1})_{p}\leq
C_mk^{-2b_mr}\sum^{k}_{v=1}v^{2b_mr-1}\|\oplus^rJ_{v,s}(f)-f\|_{p},
$$
where $0<b_{m}<1$, $b_{m}\rightarrow 1$, as $m \rightarrow \infty$.

Set
\[\sigma_{v}=v^{-2r}\|\widetilde{\Delta}^{r}\left(\oplus^rJ_{v,s}(f)\right)\|_{p},\;\;\tau_{v}=\|\oplus^rJ_{v,s}(f)-f\|_{p},\quad v
\geq 1.\]
Using Lemma~\ref{lm31}, we have
\begin{eqnarray*}
\sigma_{k}&\leq &
k^{-2r}\|\widetilde{\Delta}^{r}\left(\oplus^rJ_{k,s}(\oplus^rJ_{v,s}(f))\right)\|_{p}+k^{-2r}\|\widetilde{\Delta}^{r}\left(\oplus^rJ_{k,s}(\oplus^rJ_{v,s}(f)-f)\right)\|_{p}\\
&\leq &
2^r\left(\frac{v}{k}\right)^{2r}\left(v^{-2r}\|\widetilde{\Delta}^{r}\left(\oplus^rJ_{v,s}(f)\right)\|_{p}\right)+C\:\|\oplus^rJ_{v,s}(f)-f\|_p\\
& = &\left(\frac{\sqrt{2}\:v}{k}\right)^{2r}\sigma_{v}+C\tau_{v}.
\end{eqnarray*}
By Lemma~\ref{sec3lm3}, we have \; $$\sigma_{k}\leq
C_mk^{-2b_mr}\sum^{k}_{v=1}v^{2b_mr-1}\tau_{v}$$ for some large enough $m$.\\ That is,
\[k^{-2r}\|\widetilde{\Delta}^{r}\left(\oplus^uJ_{k,s}(f)\right)\|_{p}\leq
C_mk^{-2b_mr}\sum^{k}_{v=1}v^{2b_mr-1}\|\oplus^rJ_{v,s}(f)-f\|_{p}.\] For $k\geq 1$,
there exists a positive integer $k_0$, $\displaystyle\frac{k}{2}\leq
k_0\leq k$, such
that\[\|\oplus^rJ_{k_0,s}(f)-f\|_{p}\leq\|\oplus^rJ_{v,s}(f)-f\|_{p},\quad
\frac{k}{2}\leq v\leq k.\] Thus
\begin{eqnarray*}
K_{2r}(f,\widetilde{\Delta}, k^{-2r})_{p}&\leq &
\|\oplus^rJ_{k_0,s}(f)-f\|_{p}+k^{-2r}\|\widetilde{\Delta}^{r}\left(\oplus^rJ_{k_0,s}(f)\right)\|_{p}\\
&\leq & 2^{2r}k^{-2r}\sum_{\frac{k}{2}\leq v \leq
k}v^{2r-1}\|\oplus^rJ_{v,s}(f)-f\|_{p}\\
& & +C_m k^{-2b_mr}\sum^{k}_{v=1}v^{2b_mr-1}\|\oplus^rJ_{v,s}(f)-f\|_{p}\\
&\leq &
C_mk^{-2b_mr}\sum^{k}_{v=1}v^{2b_mr-1}\|\oplus^rJ_{v,s}(f)-f\|_{p}.
\end{eqnarray*}
From (\ref{sec2eq5}) it follows that
\begin{eqnarray*}
\omega^{2r}(f,k^{-1})_p\leq C_m
k^{-2b_mr}\sum_{v=1}^{k}v^{2b_mr-1}\|\oplus^uJ_{v,s}(f)-f\|_{p}.
\end{eqnarray*}

To finish our proof, we need the following inequalities.
\begin{eqnarray}\label{sec4eq25}
\omega^{2r}(f,k^{-1})_p & \approx & \frac{1}{k^{2r}}\max_{1\leq
v\leq k}v^{2r} \|\oplus^rJ_{v,s}(f)-f\|_p\nonumber\\ & \approx &
\frac{1}{k^{2r+\frac{1}{4}}}\max_{1\leq v\leq k}v^{2r+\frac{1}{4}}
\|\oplus^rJ_{v,s}(f)-f\|_p.
\end{eqnarray}
In the first place, we prove the former inequality of
(\ref{sec4eq25}) (explained below).
\begin{eqnarray*}
\omega^{2r}(f,k^{-1})_p &\leq &
C_1 k^{-2b_mr}\sum_{v=1}^{k}v^{2b_mr-1}\|\oplus^rJ_{v,s}(f)-f\|_{p}\\
&\leq & C_1
\left(k^{-2b_mr}\sum_{v=1}^{k}v^{-2(1-b_m)r-1}\right)\max_{1\leq
v \leq k}v^{2r}\|\oplus^rJ_{v,s}(f)-f\|_{p}\\
&\leq & C_2k^{-2r}\max_{1\leq v \leq k}v^{2r}\|\oplus^rJ_{v,s}(f)-f\|_{p}\\
&\leq & C_3k^{-2r}\max_{1\leq v \leq k}v^{2r}\omega^{2r}(f,v^{-1})_{p}\\
&\leq & C_4\left(k^{-2r}\max_{1\leq v \leq k}v^{2r}\left(\frac{k}{v}\right)^{2r}\right)\omega^{2r}(f,k^{-1})_{p}\\
&\leq & C_4\omega^{2r}(f,k^{-1})_{p},
\end{eqnarray*}
where the fourth inequality is deduced by Theorem~\ref{th41} and
the fifth is by (\ref{sec2eq6}).\\
 Thus
\[\omega^{2r}(f,k^{-1})_p\approx \frac{1}{k^{2r}}\max_{1\leq v\leq k} v^{2r}\|\oplus^rJ_{v,s}(f)-f\|_{p}.\] In the same way, we have
\begin{eqnarray*}
\omega^{2r}(f,k^{-1})_{p} & \leq &
C_1 k^{-2b_mr}\sum_{v=1}^k v^{2b_mr-1}\|\oplus^rJ_{v,s}(f)-f\|_{p}\\
& \leq & C_1 \left(k^{-2b_mr}\sum_{v=1}^k
v^{-2(1-b_m)r-\frac{1}{4}-1}\right)\max_{1\leq v\leq k}
v^{2r+\frac{1}{4}}\|\oplus^rJ_{v,s}(f)-f\|_p\\
& \leq & C_5 k^{-2r-\frac{1}{4}} \max_{1\leq v\leq k}
v^{2r+\frac{1}{4}}\|\oplus^rJ_{v,s}(f)-f\|_p\\
& \leq& C_6 k^{-2r-\frac{1}{4}}\max_{1\leq v\leq k}
v^{2r+\frac{1}{4}} \omega^{2r}(f,v^{-1})_p\\
& \leq & C_6 k^{-2r-\frac{1}{4}}\left(\max_{1\leq v\leq k}
v^{2r+\frac{1}{4}}\left(\frac{k}{v}\right)^{2r}\right)\omega^{2r}(f,k^{-1})_p\\
& \leq & C_7 \omega^{2r}(f,k^{-1})_p,
\end{eqnarray*}
that is,
\begin{eqnarray*}
\omega^{2r}(f,k^{-1})_{p}\approx \frac{1}{k^{2r+\frac{1}{4}}}\max_{v
\geq k}v^{2r+\frac{1}{4}}\|\oplus^rJ_{v,s}(f)-f\|_{p}.
\end{eqnarray*}
Therefore
\begin{eqnarray}\label{sec4eq26}
\omega^{2r}(f,k^{-1})_{p}&\approx & \frac{1}{k^{2r}}\max_{v\geq
k}v^{2r}\|\oplus^r J_{v,s}(f)-f\|_p\nonumber\\
&\approx & \frac{1}{k^{2r+\frac{1}{4}}}\max_{v \geq
k}v^{2r+\frac{1}{4}}\|\oplus^rJ_{v,s}(f)-f\|_{p}.
\end{eqnarray}

Now we can complete the proof of (\ref{sec4eq21}).
Clearly, there exists
$1\leq k_1\leq k$ such that
\begin{eqnarray*}
k_1^{2r+\frac{1}{4}}\|\oplus^rJ_{k_1,s}(f)-f\|_{p}=\max_{1\leq v
\leq k}v^{2r+\frac{1}{4}}\|\oplus^rJ_{v,s}(f)-f\|_{p}.
\end{eqnarray*}
Then it is deduced from (\ref{sec4eq26}) that
\begin{eqnarray*}
k^{-2r}k_{1}^{2r}\|\oplus^rJ_{k_{1},s}(f)-f\|_{p} &\leq &
\frac{1}{k^{2r}}\max_{1\leq v \leq
k}v^{2r}\|\oplus^rJ_{v,s}(f)-f\|_{p}\\
&\leq & C_8\frac{1}{k^{2r+\frac{1}{4}}}\max_{1 \leq v \leq
k}v^{2r+\frac{1}{4}}\|\oplus^rJ_{v,s}(f)-f\|_{p}\\
& = &
C_8k^{-2r-\frac{1}{4}}k_{1}^{2r+\frac{1}{4}}\|\oplus^rJ_{k_1,s}(f)-f\|_{p}.
\end{eqnarray*}
This implies $k_{1}\approx k$. Applying (\ref{sec4eq26}) again implies
\begin{eqnarray*}
\omega^{2r}(f,k^{-1})_p &\leq &
C_{5}\frac{1}{k^{2r+\frac{1}{4}}}\max_{1 \leq v \leq
k}v^{2r+\frac{1}{4}}\|\oplus^rJ_{v,s}(f)-f\|_{p}\\
&= &
C_{5}\frac{1}{k^{2r+\frac{1}{4}}}(k_{1}^{2r+\frac{1}{4}}\|\oplus^rJ_{k_{1},s}(f)-f\|_{p})\\
&\leq & C_{5}\max_{k_1\leq v \leq k}\|\oplus^rJ_{v,s}(f)-f\|_{p}.
\end{eqnarray*}
Noticing that $k_{1}\approx k$, we may rewrite the above inequality
as \[\omega^{2r}(f,k^{-1})_{p}\leq C\max_{v\geq
k}\|\oplus^rJ_{v,s}(f)-f\|_{p}.\]

This completes the
proof of Theorem 4.2.
\quad$\Box$

\begin{theorem}\label{th43}
$\{\oplus^rJ_{k,s}\}^{\infty}_{k=1}$ are saturated on
$L^{p}(\mathbb S^{n-1})$ with order $k^{-2r}$ and the collection of
constants is their invariant class.
\end{theorem}

{\bf Proof.} We first prove for $j=0,1,2,\dots$,
\begin{eqnarray}\label{sec4eq23}
\lim_{k\rightarrow\infty}\frac{1- {^1\xi_k(j)}}{1- {^1\xi_k(1)}} & =
& \frac{j(j+2\lambda)}{2\lambda+1}.
\end{eqnarray}
In fact, for any $0<\delta<\pi$, it follows from (\ref{eq34}) that
\begin{eqnarray*}
\int_{\delta}^{\pi}\mathscr{D}_{k,s}(\theta)\sin^{2\lambda}\theta\:d\theta
& \leq &
\int_{\delta}^{\pi}{\left(\frac{\theta}{\delta}\right)}^{\!\!3}\mathscr{D}_{k,s}(\theta)\sin^{2\lambda}\theta\:d\theta\\
& \leq &
\delta^{-3}\int_{0}^{\pi}\theta^{\:3}\mathscr{D}_{k,s}(\theta)\sin^{2\lambda}\theta\:d\theta
\:\leq\; C_{\delta,s}k^{-3}.
\end{eqnarray*}
For $v=1,2,\dots$, we have, using (\ref{eq34}) again,
\begin{eqnarray}\label{sec4eq24}
 1-{^1\xi_k(1)}
  =
 \int_0^{\pi}\mathscr{D}_{k,s}(\theta)\left(1-\frac{G_1^\lambda(\cos\theta)}{G_1^\lambda(1)}\right)\sin^{2\lambda}\theta\:d\theta
 & \approx &
 \int_0^{\pi}\mathscr{D}_{k,s}(\theta)\sin^2\frac{\theta}{2}\sin^{2\lambda}\theta\:d\theta\nonumber\\
 & \approx & k^{-2}.
\end{eqnarray}

We deduce from (\ref{eq25}) that for any $\epsilon>0$, there exists
$\delta>0$, for $0<\theta<\delta$, such that
\begin{eqnarray*}
\left|\left(1-P_j^n(\cos\theta)\right)-\frac{j(j+2\lambda)}{2\lambda+1}\left(1-P_1^n(\cos\theta)\right)\right|
& \leq & \epsilon\left(1-P_1^n(\cos\theta)\right).
\end{eqnarray*}
Then it follows that
\begin{eqnarray*}
& & \left|\left(1-{^1\xi_k(j)}\right)-\frac{j(j+2\lambda)}{2\lambda+1}\left(1-{^1\xi_k(1)}\right)\right|\\
& = &
\left|\int_{0}^{\pi}\mathscr{D}_{k,s}(\theta)\left(1-P_j^n(\cos\theta)\right)\sin^{2\lambda}\theta\:d\theta\right.\\
& &
\left.-\int_{0}^{\pi}\mathscr{D}_{k,s}(\theta)\left(1-P_1^n(\cos\theta)\right)\frac{j(j+2\lambda)}{2\lambda+1}\sin^{2\lambda}\theta\:d\theta\right|\\
& = &
\left|\int_{0}^{\pi}\mathscr{D}_{k,s}(\theta)\left(\big(1-P_j^n(\cos\theta)\big)-
\big(1-P_1^n(\cos\theta)\big)\frac{j(j+2\lambda)}{2\lambda+1}\right)\sin^{2\lambda}\theta\:d\theta\right|\\
& \leq &
\int_{0}^{\delta}\mathscr{D}_{k,s}(\theta)\:\epsilon\:\left(1-P_1^n(\cos\theta)\right)\sin^{2\lambda}\theta\:d\theta
+2 \int_{\delta}^{\pi}\mathscr{D}_{k,s}(\theta)\sin^{2\lambda}\theta\left(1+\frac{j(j+2\lambda)}{2\lambda+1}\right)\:d\theta\\
& \leq & C \epsilon\:k^{-2}+C_{\delta,s}k^{-3}.
\end{eqnarray*}
 So, (\ref{sec4eq23}) holds.
By Lemma~\ref{lm34},
\[\oplus^rJ_{k,s}(f)-f=\sum^{\infty}_{j=0}\left(1-^r\xi_{k,s}(j)\right)Y_{j}(f)\]
and for $j=1,2,\dots$,
\begin{eqnarray*}
1-{^r\xi_{k,s}}(j)=\left(\int_0^{\gamma}
\mathscr{D}_{k,s}(\theta)\left(1-P_j^n(\cos\theta)
\right)\sin^{2\lambda}\theta\:d\theta\right)^r
=\left(1-{^1\xi_{k,s}}(j)\right)^r.
\end{eqnarray*}
Combining with (\ref{sec4eq23}) and (\ref{sec4eq24}), we have,
\[\lim_{k\rightarrow
\infty}\frac{1-{^r\xi_{k,s}}(j)}{1-{^r\xi_{k,s}}(1)}=
\left(\frac{j(j+2\lambda)}{2\lambda}\right)^r\neq 0\] and
$1-{^r\xi_{k,s}}(1)\approx k^{-2r}.$ Using of Lemma~\ref{lm38}, we
finish the proof of Theorem 4.3.\quad$\Box$

We obtain the following corollary from Theorem~\ref{th41},
Theorem~\ref{th42} and Theorem~\ref{th43}.

\begin{corollary}
For positive integers $r$ and $s$, $2s\geq n$,
 $0<\alpha\leq 2r$,
 $f\in
L^{p}(\mathbb{S}^{n-1}), 1\leq p\leq \infty$, and the
sequence of  Boolean sums of Jackson
operators $\left\{\oplus^rJ_{k,s}\right\}_{k=1}^{\infty}$ given by
(\ref{sec2eq8}), the
following statements  are equivalent. \\
\begin{tabular}{ll}
(i)\quad
$\| \oplus^rJ_{k,s}(f) -
f\|_p=O(k^{-\alpha})$\quad($k\rightarrow \infty$);  \\
(ii)\quad
 $\omega^{2r}(f,\delta)_{p}=O(\delta^{\:\alpha})$\quad
($\delta\rightarrow 0$). \quad $\Box$
\end{tabular}
\end{corollary}

\end{document}